\documentclass[10pt]{amsart}
\usepackage{curvesls}
\usepackage{epic}

\theoremstyle{plain}
\newtheorem{thm}{Theorem}[section]
\newtheorem{lem}[thm]{Lemma}
\newtheorem{prop}[thm]{Proposition}
\newtheorem{cor}[thm]{Corollary}
\newtheorem{conj}[thm]{Conjecture}
\newtheorem{qn}[thm]{Question}

\theoremstyle{definition}

\newcommand{\Z}{{\bf{Z}}}
\newcommand{\Burau}{{\mathrm{Burau}}}
\newcommand{\co}{\colon\thinspace}

\title{Does the Jones Polynomial Detect the Unknot?}
\author{Stephen Bigelow}
\address{Department of Mathematics and Statistics,
         University of Melbourne,
         Parkville Victoria 3052, Australia}
\email{bigelow@unimelb.edu.au}
\date{December 2000}

\begin{document}
\begin{abstract}

We address the question: Does there exist a non-trivial knot with
a trivial Jones polynomial? To find such a knot, 
it is almost certainly sufficient 
to find a non-trivial braid on four strands in the kernel of the 
Burau representation. I will describe a computer algorithm to 
search for such a braid.

\end{abstract}
\maketitle
\section{Introduction}
\label{INTRODUCTION}

The Jones polynomial $V_K(q)$ of a knot $K$
is one of the most famous and important knot invariants.
It is not hard to construct
distinct knots with the same Jones polynomial.
However the answer to the following question remains unknown.

\begin{qn}
\label{v=1}
Does there exist a non-trivial knot $K$ with $V_K(q) \equiv 1$?
\end{qn}

This is given as Problem 1 in \cite{jones:problems}.
There have been many attempts to find such a knot.
A brute force approach was used in \cite{dasbach}
to check all knots with up to seventeen crossings.
Another approach used in \cite{apr} and \cite{jones-rolfsen}
is to start with a complicated diagram of the unknot
and apply mutations which do not alter the Jones polynomial
but may alter the knot type.

The approach described in this paper comes from the theory of braids.
Any knot $K$ can be obtained as the closure of some braid $\beta$.
The Jones polynomial of $K$ is a trace function of
the representation of $\beta$ into the Temperley-Lieb algebra.
We are therefore led to ask the following question.

\begin{qn}
Is the representation of the braid group
into the Temperley-Lieb algebra faithful?
\end{qn}

This is Problem 3 in \cite{jones:problems}.
A non-trivial braid in the kernel of the Temperley-Lieb representation
could be used to construct a knot
with Jones polynomial equal to one.
I am not aware of any proof
that the knot so obtained must be non-trivial,
but this seems unlikely to pose a problem
if a specific braid were known.
The following conjecture is therefore widely assumed to be true.

\begin{conj}
\label{tl_poly}
If the Temperley-Lieb representation of the braid group is unfaithful
then there exists a non-trivial knot with Jones polynomial equal to one.
\end{conj}

The Temperley-Lieb representation of $B_n$
appears as a summand in a larger representation
into the Hecke algebra $H(q,n)$ of type $A_{n-1}$.
We will call this latter representation the Jones representation,
although some authors use this term
for what we are calling the Temperley-Lieb representation.
The Jones representation
was used by Ocneanu in \cite{homfly}
to define a two-variable generalisation of the Jones polynomial
called the HOMFLY polynomial.
The following conjecture is also widely assumed to be true.

\begin{conj}
\label{hecke_poly}
If the Jones representation of the braid group is unfaithful
then there exists a non-trivial knot with HOMFLY polynomial equal to one.
\end{conj}

We will focus on the braid group $B_4$.
In this case
the Jones and Temperley-Lieb representations
both decompose into the Burau representation
together with some very simple representations.
Thus we have the following.

\begin{prop}
\label{equifaithful}
The following are equivalent:
\begin{itemize}
\item the Jones representation of $B_4$ is faithful,
\item the Temperley-Lieb representation of $B_4$ is faithful, and
\item the Burau representation of $B_4$ is faithful.
\end{itemize}
\end{prop}

We are therefore led to ask the following question.

\begin{qn}
Is the Burau representation of $B_4$ faithful?
\end{qn}

A negative answer would almost certainly lead to 
a non-trivial knot whose HOMFLY polynomial is equal to one.
As far as I know, 
a positive answer would have no such dramatic consequences
other than finally determining
for which values of $n$ the Burau representation of $B_n$ is faithful.
Krammer \cite{krammer} has already shown that $B_4$ is be linear.

The Burau representation of $B_n$
is known to be faithful for $n \le 3$ \cite{birman}
and unfaithful for $n \ge 5$ \cite{bigelow:b5}.
The case $n=4$ seems to lie
very close to the border between faithfulness and unfaithfulness.

The main aim of this paper
is to propose a computer search for 
a non-trivial braid in the kernel of the Burau representation of $B_4$.
This might seem overly ambitious.
After all, it amounts to a search for a very special case of
a non-trivial knot whose HOMFLY polynomial is equal to one 
(assuming Conjecture \ref{hecke_poly}).
Many people have tried and failed
to find a non-trivial knot
whose weaker Jones polynomial is equal to one.
However there is some reason for optimism.
A knot constructed by the methods of this paper
would have thousands of crossings.
Thus we are searching in relatively unexplored territory
which might contain unexpected treasures.
This is probably enough to justify
the expenditure of some computer time,
but perhaps not too much human time or brain power.

\section{The Burau Representation}

We now define the braid groups $B_n$
and the Burau representation.

Let $D$ be a disk.
Let $p_1,\dots,p_n$ be distinct points in the interior of $D$.
We call these ``puncture points''.
Let $D_n = D \setminus \{p_1,\dots,p_n\}$.
Let $d_0$ be a basepoint on $\partial D_n$.
For concreteness,
take $D$ to be the unit disk in the complex plane centred at the origin,
take $p_1,\dots,p_n$ to be real numbers
satisfying $-1 < p_1 < \dots < p_n < 1$,
and take $d_0$ to be $-i$.

The braid group $B_n$ is defined to be
the group of homeomorphisms from $D_n$ to itself
which act as the identity on $\partial D_n$,
taken up to isotopy relative to $\partial D_n$.
It is generated by $\sigma_1,\dots,\sigma_{n-1}$,
where $\sigma_i$ exchanges $p_i$ and $p_{i+1}$
by a counterclockwise half twist.

The fundamental group $\pi_1(D_n,d_0)$ 
is a free group with basis $x_1,\dots,x_n$,
where $x_i$ is a loop based at $d_0$
which passes counterclockwise around $p_i$ 
and no other puncture points.
Let $\phi \co \pi_1(D_n,d_0) \rightarrow \langle q \rangle$
be the homomorphism given by $\phi(x_i) = q$.
Let $\tilde{D}_n$ be the covering space
corresponding to the subgroup $\ker(\phi)$ of $\pi_1(D_n)$.
Fix a point $\tilde{d}_0$ in the fibre over $d_0$.

A more concrete description of $\tilde{D}_n$
can be given as follows.
Make a bi-infinite stack of $\Z$ copies of $D_n$.
On each copy, make a series of vertical cuts
connecting each of the puncture points $p_i$ to the boundary.
Glue the left-hand side of each cut
to the right-hand side of the corresponding cut
on the copy of $D_n$ one level lower.

The group of covering transformations of $\tilde{D}_n$
is $\langle q \rangle$.
The $\Z$-module $H_1(\tilde{D}_n)$
can be considered as a $\Z[q^{\pm 1}]$-module,
where multiplication by $q$
is the induced action of the covering transformation $q$.
Thought of in this way,
$H_1(\tilde{D}_n)$ turns out to be
a free $\Z[q^{\pm 1}]$-module of rank $n-1$.

The Burau representation
is the induced action of $B_n$
by $\Z[q^{\pm 1}]$-module homomorphisms
on $H_1(\tilde{D}_n)$.
We make this more precise as follows.
Let $\beta \co D_n \rightarrow D_n$ be a homeomorphism
representing a braid $[\beta]$ in $B_n$.
The induced action of $\beta$ on $\pi_1(D_n)$
satisfies $\phi \beta = \phi$.
It follows by some basic algebraic topology
that there exists a unique lift $\tilde{\beta}$
which makes the following diagram commute.
$$
\begin{array}{ccc}
(\tilde{D}_n,\tilde{d}_0)
  & \stackrel{\tilde{\beta}}{\rightarrow} 
  & (\tilde{D}_n,\tilde{d}_0) \\
\downarrow  &                               & \downarrow \\
(D_n,d_0)   & \stackrel{\beta}{\rightarrow} & (D_n,d_0)
\end{array}
$$
Furthermore, $\tilde{\beta}$ commutes with
the action of $q$ on $\tilde{D}_n$ by a covering transformation.
Thus $\tilde{\beta}$ induces a $\Z[q^{\pm 1}]$-module homomorphism
$$
\tilde{\beta}_* \co H_1(\tilde{D}_n) \rightarrow H_1(\tilde{D}_n).
$$
The Burau representation is the map
$$\Burau([\beta]) = \tilde{\beta}_*.$$

For example,
using an appropriate choice of basis for $H_1(\tilde{D}_4)$,
the Burau representation of $B_4$ is given by
\begin{eqnarray*}
\sigma_1 &\mapsto&
\left( \begin{array}{ccc} -q&q&0\\0&1&0\\0&0&1 \end{array} \right),
\\
\sigma_2 &\mapsto&
\left( \begin{array}{ccc} 1&0&0\\1&-q&q\\0&0&1 \end{array} \right),
\\
\sigma_3 &\mapsto&
\left( \begin{array}{ccc} 1&0&0\\0&1&0\\0&1&-q \end{array} \right).
\end{eqnarray*}

\section{The Jones and HOMFLY polynomials}
\label{jones}

In this section we make the connection between
the Burau representation of $B_4$
and the Jones and HOMFLY polynomials of a knot.
We start by proving Proposition \ref{equifaithful},
which states that in the case of $B_4$
the Jones, Temperley-Lieb, and Burau representations
are either all faithful or all unfaithful.

\begin{proof}[Proof of Proposition \ref{equifaithful}]
We will not define the Temperley-Lieb or Jones representations,
but we will use some of their basic properties,
all of which can be found in \cite{jones}.

The Jones representation of $B_n$
can be decomposed into irreducible summands,
one corresponding to each Young diagram with $n$ boxes.
Let $V_\lambda$ denote the representation
corresponding to the Young diagram $\lambda$.
The Temperley-Lieb representation
is the sum of those $V_\lambda$ for which $\lambda$ has one or two rows.

The Young diagrams with $4$ boxes are
$(4)$, $(3,1)$, $(2,2)$, $(2,1,1)$ and $(1,1,1,1)$.
The $V_\lambda$ which lie in the Temperley-Lieb representation
are as follows.
\begin{itemize}
\item $V_{(4)}$ is one-dimensional.
\item $V_{(3,1)}$ is the Burau representation.
\item $V_{(2,2)}$
   can be defined by composing the Burau representation of $B_3$
   with the map from $B_4$ to $B_3$ given by
   $\sigma_1 \mapsto \sigma_1$,
   $\sigma_2 \mapsto \sigma_2$, and
   $\sigma_3 \mapsto \sigma_1$.
\end{itemize}

The Young diagram $(2,1,1)$
is a reflection of $(3,1)$.
Reflection of the Young diagram has the effect of substituting
$\sigma_i \mapsto -q\sigma_i^{-1}$ in the corresponding representation.
It is shown in \cite{long-paton}
that the kernel of the Burau representation
is invariant under this substitution.
Thus the kernel of $V_{(2,1,1)}$
is the same as that of the Burau representation.
Finally, $V_{(1,1,1,1)}$ is one-dimensional.

If the Burau representation of $B_4$ is faithful
then so are the Temperley-Lieb and Jones representations.
Conversely, suppose $\beta$ is a non-trivial braid
in the kernel of the Burau representation of $B_4$.
Consider the commutator $[(\sigma_1 \sigma_2)^3,\beta]$.
This lies in the kernel of $V_{(3,1)}$, and hence $V_{(2,1,1)}$.
Since this is a commutator,
it lies in the kernel of any one-dimensional representation.
Since $(\sigma_1 \sigma_2)^3$ is central in $B_3$,
it also lies in the kernel of the representation corresponding to $(2,2)$.
Thus it lies in the kernel
of the Jones and Temperley-Lieb representations.

It remains to show that $[(\sigma_1 \sigma_2)^3,\beta]$ 
must be non-trivial.
This is not difficult,
but would take us too far afield.
We therefore omit this part of the proof.
\end{proof}

Suppose $\beta$ lies in the kernel of
the Temperley-Lieb representation.
Then the closures of the braids
$$\beta \sigma_1 \sigma_2 \dots \sigma_{n-1}$$
and
$$\sigma_1 \sigma_2 \dots \sigma_{n-1}.$$
have the same Jones polynomials.
The closure of the braid $\sigma_1 \sigma_2 \dots \sigma_{n-1}$ is the unknot.
Thus the closure of $\beta \sigma_1 \sigma_2 \dots \sigma_{n-1}$
has Jones polynomial equal to one.
If we could be sure that this was a non-trivial knot
then Conjecture \ref{tl_poly} would be proved.
If it is the unknot then all is not lost, since
we could use any power $\beta^k$ in place of $\beta$.
Thus the following conjecture
implies Conjecture \ref{tl_poly}.

\begin{conj}
Let $\beta$ be a non-trivial braid in $B_n$.
There exists some integer $k$
such that the closure of $\beta^k\sigma_1\sigma_2 \dots \sigma_{n-1}$
is a non-trivial knot.
\end{conj}

As well as powers of $\beta$,
we also have products of conjugates of $\beta$ at our disposal.
And in place of $\sigma_1\sigma_2\dots\sigma_{n-1}$
we could use any braid whose closure is the unknot.
Thus we can weaken the above conjecture to the following.

\begin{conj}
\label{true}
Let $H$ be a non-trivial normal subgroup of $B_n$.
Then there exists $\beta_1 \in H$ and $\beta_2 \in B_n$
such that the closure of $\beta_2$ is the unknot
but the closure of $\beta_1\beta_2$ is a non-trivial knot.
\end{conj}

The above discussion applies equally well to
the Jones representation and the HOMFLY polynomial.
Thus Conjecture \ref{true} also implies Conjecture \ref{hecke_poly}.

A counterexample to Conjecture \ref{true}
would be truly astonishing,
implying an unprecedented correlation
between the algebraic structure of $H$
and the geometric structure of the knots constructed.
However it might be quite difficult to prove this ``obvious'' conjecture.
This problem would probably be easily overcome
in the case of a specific non-trivial braid
in the kernel of the Burau representation of $B_4$.

\section{The case $n=3$}

The aim of this section is to prove the following.

\begin{thm}
\label{burau3}
The Burau representation of $B_3$ is faithful.
\end{thm}

This is a well-known result
and has been proved in many different ways
(see, for example, \cite{birman}).
The proof given here is a warm-up for the ideas that will be used later.

A {\em fork}
is an embedded tree $F$ in $D$
with four vertices $d_0$, $p_i$, $p_j$ and $z$
such that 
\begin{itemize}
\item $F$ meets the puncture points only at $p_i$ and $p_j$,
\item $F$ meets the $\partial D_n$ only at $d_0$, and
\item all three edges of $F$ have $z$ as a vertex.
\end{itemize}
The edge of $F$ which contains $d_0$ is called
the {\em handle} of $F$.
The union of the other two edges forms a single edge 
which we call the {\em tine edge} of $F$
and denote by $T(F)$.
Orient $T(F)$ so that the handle of $F$ lies to the right of $T(F)$.

A {\em noodle}
is an embedded oriented edge $N$ in $D_n$ such that
\begin{itemize}
\item $N$ goes from $d_0$ to another point on $\partial D_n$,
\item $N$ meets $\partial D_n$ only at its endpoints, and
\item a component of $D_n \setminus N$ contains precisely one puncture point.
\end{itemize}
This last requirement
was not included in the definition given in \cite{bigelow:bgal}.
Without it, Theorem \ref{faithful} is not true,
as far as I know.

Let $F$ be a fork and let $N$ be a noodle.
We define a pairing $\langle N,F \rangle$ in $\Z[q^{\pm 1}]$ as follows.
If necessary, apply a preliminary isotopy of $F$
so that $T(F)$ intersects $N$ transversely.
Let $z_1,\dots,z_k$ denote 
the points of intersection between $T(F)$ and $N$
(in no particular order).
For each $i=1,\dots,k$,
let $\gamma_i$ be the arc in $D_n$
which goes from $d_0$ to $z_i$ along $F$,
then back to $d_0$ along $N$.
Let $a_i$ be the integer such that $\phi(\gamma_i) = q^{a_i}$.
In other words, $a_i$ is the sum of
the winding numbers of $\gamma_i$
around each of the puncture points $p_j$.
Let $\epsilon_i$ be the sign of the intersection
between $N$ and $F$ at $z_i$.
Let
\begin{equation}
\label{pair1}
\langle N,F \rangle = \sum_{i=1}^k \epsilon_i q^{a_i}.
\end{equation}

We should really check that this is independent of 
our choice of preliminary isotopy of $F$.
This is easy enough to prove directly.
It is also a special case of the following lemma.

\begin{lem}[The Basic Lemma]
Let $\beta \co D_n \rightarrow D_n$
represent an element of the kernel of the Burau representation.
Then $\langle N,F \rangle = \langle N,\beta(F) \rangle$
for any noodle $N$ and fork $F$.
\end{lem}

\begin{proof}
We can assume that the tine edges of $F$ and $\beta(F)$
both intersect $N$ transversely.

Let $\tilde{F}$ be the lift of $F$ to $\tilde{D}_n$
which contains $\tilde{d}_0$.
Let $\tilde{T}(F)$ be the corresponding lift of $T(F)$.
Then $\tilde{T}(F)$ intersects $q^a \tilde{N}$
transversely for any $a \in \Z$.
Let $(q^a \tilde{N}, \tilde{T}(F))$
denote the algebraic intersection number of these two arcs.
Then the following definition of $\langle N,F \rangle$
is equivalent to Equation (\ref{pair1}).
\begin{equation}
\label{pair2}
\langle N,F \rangle = \sum_{a \in \Z} (q^a \tilde{N}, \tilde{T}(F))q^a.
\end{equation}

Suppose $T(F)$ goes from $p_i$ to $p_j$.
Let $\nu(p_i)$ and $\nu(p_j)$
be disjoint small regular neighbourhoods of $p_i$ and $p_j$ respectively.
Let $\gamma$ be a subarc of $T(F)$
which starts in $\nu(p_i)$ and ends in $\nu(p_j)$.
Let $\delta_i$ be a loop in $\nu(p_i)$ based at $\gamma(0)$
which passes counterclockwise around $p_i$.
Similarly, let $\delta_j$ be a loop in $\nu(p_j)$ based at $\gamma(1)$
which passes counterclockwise around $p_j$.
Let $T_2(F)$ be the ``figure eight''
$$T_2(F) = \gamma \delta_j \gamma^{-1} \delta_i^{-1}.$$
Let $\tilde{T}_2(F)$ be the lift of $T_2(F)$
which is equal to $(1-q)\tilde{T}(F)$ 
outside a small neighbourhood of the puncture points.
Then the following definition of $\langle N,F \rangle$
is equivalent to Equation (\ref{pair2}).
\begin{equation}
\label{pair3}
\langle N,F \rangle = \frac{1}{1-q}
\sum_{a \in \Z} (q^a \tilde{N}, \tilde{T}_2(F)) q^a.
\end{equation}

Note that $\tilde{T}_2(F)$ is a closed loop in $\tilde{D}_n$.
Since $\beta$ is in the kernel of the Burau representation,
the loops $\tilde{T}_2(F)$ and $\tilde{T}_2(\beta(F))$ 
represent the same element of $H_1(\tilde{D}_n)$.
They therefore have the same algebraic intersection number
with any lift $q^a \tilde{N}$ of $N$.
Thus Equation \ref{pair3} will give the same result
for $\langle N,\beta(F) \rangle$ as for $\langle N,F \rangle$.
\end{proof}

We now use the assumption that $n=3$.

\begin{lem}[The Key Lemma]
In the case $n=3$,
$\langle N,F \rangle = 0$
if and only if
$T(F)$ is isotopic to an arc which is disjoint from $N$.
\end{lem}

\begin{proof}
Apply an isotopy to $F$ so that 
$T(F)$ intersects $N$ at a minimum number of points,
which we denote $z_1,\dots,z_k$ (in no particular order).
Recall the definition given in Equation (\ref{pair1}).
$$\langle N,F \rangle = \sum_{i=1}^k \epsilon_i q^{a_i}.$$
If $k=0$ then clearly $\langle N,F \rangle = 0$.
We now assume that $k>0$ and prove that $\langle N,F \rangle \neq 0$.

By applying a homeomorphism to our picture,
we can take $N$ to be
a horizontal straight line through $D_3$
with two puncture points above it
and one puncture point below it.
(The noodle has been pulled straight
and the fork is twisted!)
Let $D_n^+$ and $D_n^-$ be
the upper and components of $D_n \setminus N$ respectively.
Relabel the puncture points so that
$D_n^+$ contains $p_1$ and $p_2$
and $D_n^-$ contains $p_3$.

Consider the intersection of $T(F)$ with $D_n^-$.
This consists of a disjoint collection of arcs
which have both endpoints on $N$,
and possibly one arc with an endpoint on $p_3$.
An arc in $T(F) \cap D_n^-$ which has both endpoints on $N$ 
must enclose $p_3$,
since otherwise it could be slid off $N$
to reduce the number of points of intersection between $T(F)$ and $N$.
Thus $T(F) \cap D_n^-$ must consist of
a collection of parallel arcs enclosing $p_3$,
and possibly one arc with an endpoint on $p_3$.

Similarly, each of the arcs in $T(F) \cap D_n^+$
either enclose one of the puncture points $p_1$ or $p_2$,
or have an endpoint on one of $p_1$ or $p_2$.
There can be no arc in $T(F) \cap D_n^+$
which encloses both $p_1$ and $p_2$,
since the outermost such arc
together with the outermost arc in $T(F) \cap D_n^-$
would form a closed loop.

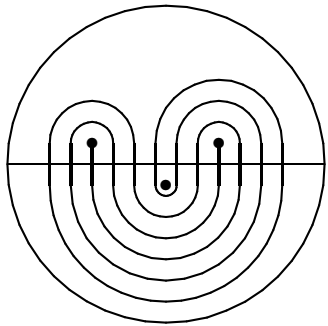
\begin{figure}
\centering
\begin{picture}(120,120)
\thicklines
\put(60,60){\bigcircle{120}}

\thinlines
\put(0,60){\line(1,0){120}}

\thicklines

\multiput(16,52)(8,0){12}{\line(0,1){16}}

\put(32,68){\arc(8,0){180}}
\put(32,68){\arc(16,0){180}}

\put(80,68){\arc(8,0){180}}
\put(80,68){\arc(16,0){180}}
\put(80,68){\arc(24,0){180}}

\put(60,52){\arc(-4,0){180}}
\put(60,52){\arc(-12,0){180}}
\put(60,52){\arc(-20,0){180}}
\put(60,52){\arc(-28,0){180}}
\put(60,52){\arc(-36,0){180}}
\put(60,52){\arc(-44,0){180}}

\put(32,68){\circle*{4}}
\put(80,68){\circle*{4}}
\put(60,52){\circle*{4}}
\end{picture}
\caption{A tine edge and a noodle in $D_3$.}
\label{d3}
\end{figure}
An example of a noodle and a tine edge in $D_3$
is shown in Figure \ref{d3}.
We have omitted the handle of the fork,
which plays no role in our argument.

Let $z_i$ and $z_j$ be two points of intersection
between $T(F)$ and $N$
which are joined by an arc in
$T(F) \cap D_n^+$ or $T(F) \cap D_n^-$.
This arc, together with a subarc of $N$, encloses one puncture point.
Thus
$$a_j = a_i \pm 1.$$
Also, $T(F)$ intersects $N$ with opposite signs at $z_i$ and $z_j$, so
$$\epsilon_j = -\epsilon_i.$$
Thus
$$\epsilon_j(-1)^{a_j} = \epsilon_i(-1)^{a_i}.$$

Proceeding along $T(F)$, we conclude that
the values of $\epsilon_i(-1)^{a_i}$ are 
the same for all $i=1,\dots,k$.
Thus $\langle N,F \rangle$ evaluated at $q=-1$
is equal to $\pm k$.
Thus $\langle N,F \rangle$ is not equal to zero.
\end{proof}

We are now ready to prove that 
the Burau representation of $B_3$ is faithful.

\begin{proof}[Proof of Theorem \ref{burau3}]
Let $\beta \co D_3 \rightarrow D_3$ be a homeomorphism
which represents an element of the kernel of the Burau representation.
We will show that $\beta$ is 
isotopic relative to $\partial D_n$ to the identity map,
and so represents the trivial braid.

Let $N$ be a noodle.
As before,
take $N$ to be a horizontal line through $D_n$
such that the puncture points $p_1$ and $p_2$
lie above $N$
and $p_3$ lies below $N$.
Let $F$ be a fork such that $T(F)$ is a straight line
from $p_1$ to $p_2$ which does not intersect $N$.
Then $\langle N,F \rangle = 0$.
By the Basic Lemma, $\langle N,\beta(F) \rangle = 0$.
By the Key Lemma, $\beta(T(F))$ is isotopic
to an arc which is disjoint from $N$.
By applying an isotopy to $\beta$ relative to $\partial D_n$,
we can assume that $\beta(T(F)) = T(F)$.

By a similar argument using different noodles,
we can assume that $\beta$ fixes
the triangle with vertices $p_1$, $p_2$ and $p_3$.
Thus $\beta$ must be some power of
$\Delta$, the Dehn twist about a curve parallel to $\partial D_n$.
It is easy to show that the Burau representation of $\Delta$
is the scalar matrix $q^3 I$.
Thus the only power of $\Delta$
which lies in the kernel of the Burau representation 
is the trivial braid.
\end{proof}

\section{The case $n=4$}

We now address the question
of whether the Burau representation of $B_4$ is faithful.
If the Key Lemma holds for the case $n=4$
then the same argument used for $B_3$
can be used to show that the Burau representation of $B_4$ is faithful.
The converse is also true:
if the Key Lemma is false for a given $n$
then the Burau representation of $B_n$ is unfaithful.
In other words, the following theorem holds.

\begin{thm}
\label{faithful}
The following are equivalent:
\begin{itemize}
\item the Burau representation of $B_n$ is faithful,
\item if $N$ and $F$ are any noodle and fork in $D_n$
  such that $\langle N,F \rangle = 0$
  then $T(F)$ is isotopic to an arc which is disjoint from $N$.
\end{itemize}
\end{thm}

A proof can be found in \cite{bigelow:b5},
although the terminology of noodles and forks is not used.
The proof of one direction is much the same as
our proof that the Burau representation of $B_3$ is faithful.
The proof of the other direction is constructive.
Suppose $\langle N,F \rangle = 0$
but $T(F)$ is not isotopic to an arc which is disjoint from $N$.
Let $\gamma_1$ be a simple closed curve
which is parallel to the boundary of
the component of $D_n \setminus N$
containing all but one puncture point.
Let $\gamma_2$ be the boundary of a regular neighbourhood of $T(F)$.
It is shown that 
the commutator of the Dehn twists about $\gamma_1$ and $\gamma_2$
is a non-trivial braid in the kernel of the Burau representation of $B_n$.

We now define a {\em standard form}
for a noodle $N$ and tine edge $T(F)$,
similar to the one used in the proof of the Key Lemma.
Let $N$ be a horizontal straight line through $D_4$
with $p_1$, $p_2$ and $p_3$ above it, and $p_4$ below it.
Let $D_4^+$ and $D_4^-$ be 
the upper and lower halves of $D_4 \setminus N$, respectively.
Then $D_4^- \cap T(F)$ is
a collection of disjoint arcs which enclose $p_4$,
and possibly one arc with an endpoint on $p_4$.
Each arc in $D_4^+ \cap T(F)$ either
\begin{itemize}
\item encloses one of $p_1$, $p_2$ or $p_3$,
\item encloses $p_1$ and $p_2$,
      the two leftmost puncture points in $D_4^+$, or
\item has an endpoint on a puncture point.
\end{itemize}
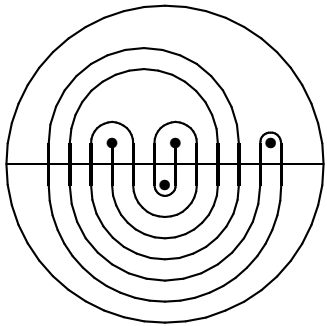
\begin{figure}
\centering
\begin{picture}(120,120)
\thicklines
\put(60,60){\bigcircle{120}}

\thinlines
\put(0,60){\line(1,0){120}}

\thicklines

\multiput(16,52)(8,0){12}{\line(0,1){16}}

\put(40,68){\arc(8,0){180}}
\put(64,68){\arc(8,0){180}}
\put(52,68){\arc(36,0){180}}
\put(52,68){\arc(28,0){180}}

\put(100,68){\arc(4,0){180}}

\put(60,52){\arc(-4,0){180}}
\put(60,52){\arc(-12,0){180}}
\put(60,52){\arc(-20,0){180}}
\put(60,52){\arc(-28,0){180}}
\put(60,52){\arc(-36,0){180}}
\put(60,52){\arc(-44,0){180}}

\put(40,68){\circle*{4}}
\put(64,68){\circle*{4}}
\put(100,68){\circle*{4}}
\put(60,52){\circle*{4}}
\end{picture}
\caption{A tine edge and a noodle in $D_4$.}
\label{d4}
\end{figure}
Figure \ref{d4} shows
an example of a noodle and a tine edge in standard form in $D_4$.

Any noodle $N$ and tine edge $T(F)$ can be put into standard form
by first isotoping $T(F)$ so as to intersect $N$
at a minimum possible number of points,
and then applying some homeomorphism to the entire picture.
The homeomorphism might need to be orientation-reversing.
This would have the effect of
substituting $q^{-1}$ for $q$ in $\langle N,F \rangle$,
so would not affect whether $\langle N,F \rangle$ is zero.

The simple parity argument used to prove the Key Lemma in $D_3$
will not work for $D_4$ because of
the existence of arcs enclosing two puncture points.
In fact, in $D_4$ there can be some cancellation
in the calculation of $\langle N,F \rangle$,
whereas our argument showed that this cannot happen in $D_3$.
We might attempt a more sophisticated argument
which shows that there cannot be {\em complete} cancellation.
Unfortunately, none of the obvious approaches seem to work.
For example, it is possible to have complete cancellation of
all of the highest and lowest powers of $q$
that occur in the calculation of $\langle N,F \rangle$.

Conversely, we could attempt a computer search
to find a counterexample to the Key Lemma for $n=4$,
and hence a non-trivial braid in the kernel of
the Burau representation of $B_4$.
This approach has worked for $B_5$ \cite{bigelow:b5}.

A tine edge $T(F)$ in standard form
is determined up to isotopy by the following:
\begin{itemize}
\item
four non-negative integers
specifying the number of arcs in $T(F) \cap D_n^+$
of each of the four possible types, and
\item
which of the puncture points are endpoints of $T(F)$.
\end{itemize}
The handle of $F$ can be ignored
because it has no effect on $\langle N,F \rangle$
up to sign and multiplication by a power of $q$.

By some of the basic theory of curves on surfaces,
if $T(F)$ is in standard form
and intersects $N$
then it is not isotopic to an arc which is disjoint from $N$.
Given data defining $T(F)$,
it is easy to compute $\langle N,F \rangle$
up to sign and multiplication by a power of $q$.
We can thus embark upon an exhaustive open-ended search
for a tine edge $T(F)$ in standard form
which intersects $N$ but gives $\langle N,F \rangle = 0$.
We now discuss issues of speed.

The polynomial $\langle N,F \rangle$ can be stored as
an array of integers.
Working with this array takes a significant amount of computer time.
There is a simple trick which can be used to eliminate this problem.
Let $M$ be a large integer.
Consider a map
$$\Z[q^{\pm 1}] \rightarrow \Z/M\Z$$
sending $q$ to some unit in $\Z/M\Z$.
Instead of computing $\langle N,F \rangle$
we can compute its image in $\Z/M\Z$.
This allows us to work with a single integer instead of an array.
There will be some ``false alarms''
for which $\langle N,F \rangle$ is non-zero
but its image in $\Z/M\Z$ is zero.
However these are infrequent and easily checked separately.

This trick speeds up the search considerably.
I have used it to check all forks for which $T(F)$ intersects $N$ at
up to $2000$ points.
By comparison, the example in $D_5$
consists of a noodle and a tine edge
which intersect at $100$ points.

There are some possibilities for further improvements in the algorithm.
Perhaps the simplest way to speed up the search
is to increase the number of searchers.
I would like to take this opportunity
to advertise my webpage
\begin{quote}
http://www.ms.unimelb.edu.au/$\sim$bigelow
\end{quote}
where, at the time of writing,
it is possible to donate computer time to this
noble and possibly futile search.

\section{Specialising $q$}

We conclude this paper with an aside
concerning the ``false alarms'' mentioned in the previous section.
Recall that a false alarm occurs when
$\langle N,F \rangle$ is non-zero
but maps to zero in $\Z/M\Z$
when $q$ is assigned some unit $q_0$.
Usually this is not very interesting,
since $M$ was fairly arbitrary.
But some false alarms occur
when the integer $q_0$ is a root of $\langle N,F \rangle$.
At first I thought that these more interesting false alarms
should give rise to a non-trivial element of the kernel
of the specialisation of the Burau representation to $q=q_0$.
However it turns out that the correct theorem is as follows.

\begin{thm}
\label{special}
Let $q_0$ be a complex number
which is not zero or a root of unity.
The following are equivalent:
\begin{itemize}
\item the Burau representation of $B_n$ is faithful
  when $q$ is specialised to $q_0$,
\item if $N$ and $F$ are any noodle and fork in $D_n$
  such that both $q_0$ and $1/q_0$ are roots of $\langle N,F \rangle$
  then $T(F)$ is isotopic to an arc which is disjoint from $N$.
\end{itemize}
\end{thm}

A computer search took about half a minute to find the following.

\begin{cor}
The Burau representation of $B_4$
is not faithful at $q=2$.
\end{cor}

\begin{proof}
\begin{figure}
\centering
\begin{picture}(120,120)
\thicklines
\put(60,60){\bigcircle{120}}

\thinlines
\put(0,60){\line(1,0){120}}

\thicklines

\multiput(20,60)(10,0){3}{\line(0,1){10}}
\multiput(50,60)(20,0){2}{\line(0,1){10}}
\multiput(80,60)(10,0){3}{\line(0,1){10}}
\put(20,50){\line(0,1){10}}
\put(100,50){\line(0,1){10}}

\multiput(30,70)(30,0){3}{\arc(10,0){180}}
\put(30,82){\makebox(0,0)[b]{$24$}}
\put(60,82){\makebox(0,0)[b]{$18$}}
\put(90,82){\makebox(0,0)[b]{$11$}}

\put(60,50){\arc(-40,0){180}}
\put(60,9){\makebox(0,0)[t]{$54$}}

\multiput(30,70)(30,0){3}{\circle*{4}}
\put(60,50){\circle*{4}}
\end{picture}
\caption{A tine edge and a noodle in $D_4$.}
\label{q=2}
\end{figure}
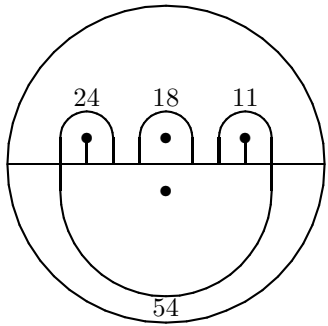
Let $T(F)$ be the tine edge in standard form
as shown schematically in Figure \ref{q=2}.
The endpoints of $T(F)$ at $p_1$ and $p_3$ are shown.
Segments of $T(F)$ are labelled with numbers
to indicate the number of parallel copies required.

A laborious computation
or a short computer program can be used to check that
$$\langle N,F \rangle = -(q - 1) (q - 2) (2 q - 1) (q^2- q + 1) (q^2+ 1),$$
up to multiplication by a power of $q$.
Both $2$ and $1/2$ are roots of this polynomial.
\end{proof}

We can construct a specific non-trivial braid $\beta$
in the kernel of the Burau representation of $B_4$ at $q=2$.
To make things more readable, let
$a=\sigma_1$, $b=\sigma_2$, and $c=\sigma_3$.
Then
$$[(ba)^3, \psi^{-1} b \psi],$$
where
\[
\psi = 
a^{-3} b^{-2} c^{-1} b c^4 b^{-1} c b a b c^2 b a^{-1} b^{-1} c^{-2}.
\]
Note, this uses
the convention that braids compose from right to left.

The noodle and fork shown in Figure \ref{q=2}
are the simplest possible example
in the sense that they have 
the fewest points of intersection.
They also have the curious property that
none of the subarcs of $T(F)$ above $N$ enclose two puncture points,
so there is no cancellation in the calculation of $\langle N,F \rangle$.
I can think of no explanation for this.

The Burau representation of $B_4$
is also unfaithful at $1/2$ and at any root of unity.
Despite hundreds of hours of computer time
I know of no other values at which it is unfaithful,
and certainly none at which it is faithful.

This is to be contrasted with the situation for $B_3$,
where we have the following.

\begin{lem}
\label{monic}
Let $N$ and $F$ be a noodle and a fork
such that $T(F)$ is not isotopic to an arc which is disjoint from $N$.
Then the highest and lowest powers of $q$
in the polynomial $\langle N,F \rangle$
both occur with coefficient $\pm 1$.
\end{lem}

\begin{cor}
If the Burau representation of $B_3$ is unfaithful at $q=q_0$.
then both $q_0$ and $1/q_0$ are roots of a monic polynomial.
In particular, the Burau representation of $B_3$
is faithful at any rational number other than $0$ or $\pm 1$.
\end{cor}

\begin{proof}[Proof of Lemma \ref{monic}]
Put $N$ and $F$ in the standard form as in Figure \ref{d3}.
Thus $N$ is a horizontal straight line
with $p_1$ and $p_2$ above it
and $p_3$ below it.
Assume that $d_0$ is the left endpoint of $N$.
We show that the lowest power of $q$ in $\langle N,F \rangle$
occurs with coefficient $\pm 1$.
The highest power of $q$
and the case where $d_0$ is the right endpoint of $N$
are handled similarly.

Let $z_1,\dots,z_k$ be the points of intersection between $N$ and $T(F)$.
Recall Equation (\ref{pair1}), which states that
$$\langle N,F \rangle = \sum_{i=1}^k \epsilon_i q^{a_i}.$$
Let $z_i$ be such that $a_i$ is minimal.
We will show that there is only one such $z_i$.
We proceed by induction $k$. 
The case $k=1$ is trivial, so assume $k>1$.

If $z_i$ were to the right of $p_3$
then there would be a subarc of $T(F)$
going from $z_i$ around $p_3$
in the clockwise (negative) sense
to intersect $N$ at a point $z_j$.
Then $a_j = a_i-1$, which contradicts the minimality of $a_i$.
Thus $z_i$ must lie to the left of $p_3$.

Let $P$ be a vertical line from the top of the disk
to a point on $N$ between the puncture points $p_1$ and $p_2$
such that $P$ does not intersect $T(F)$.
If $z_i$ were to the left of $p_3$ but to the right of $P$
then there would be a subarc of $T(F)$
going from $z_i$ around $p_2$
in the clockwise sense,
once again contradicting the minimality of $a_i$.
Thus $z_i$ lies to the left of $P$.

Let $N'$ be the union of $P$
with the portion of $N$ which lies to the left of $P$.
This is a noodle
which intersects $T(F)$ at fewer than $k$ points.
The pairing $\langle N',F \rangle$
is the sum of those monomials $\epsilon_j q^{a_j}$
for which $z_j$ lies to the left of $P$.
Thus $z_i$ is such that $a_i$ is minimal
in the calculation of $\langle N',F \rangle$.
By the induction hypothesis,
there is only one such $z_i$,
so we are done.
\end{proof}

\end{document}